\documentclass{amsart}
\usepackage{amssymb}
\usepackage{graphicx}
\usepackage[cmtip,all]{xy}

\DeclareFontFamily{U}{wncy}{}
\DeclareFontShape{U}{wncy}{m}{n}{<->wncyr10}{}
\DeclareSymbolFont{mcy}{U}{wncy}{m}{n}
\DeclareMathSymbol{\Sh}{\mathord}{mcy}{"58} 
\DeclareMathOperator{\sgn}{sign}

\newtheorem{theorem}{Theorem}[section]
\newtheorem{lemma}[theorem]{Lemma}
\newtheorem{proposition}[theorem]{Proposition}

\newtheorem{corollary}[theorem]{Corollary}

\theoremstyle{definition}

\newtheorem*{theorem*}{Theorem}
\theoremstyle{remark}

\numberwithin{equation}{section}

\begin{document}

\title{Class number divisibility for imaginary quadratic fields}


\author{Olivia Beckwith}
\address{}
\curraddr{}
\email{}
\thanks{}

\date{}

\maketitle 
\begin{abstract}
In this note we revisit classic work of Soundararajan on class groups of imaginary quadratic fields. Let $A,B,g \ge 3$ be positive integers such that $\gcd(A,B)$ is square-free. We refine Soundararajan's result to show that if $4 \nmid g$ or if $A$ and $B$ satisfy certain conditions, then the number of negative square-free $D \equiv A \pmod{B}$ down to $-X$ such that the ideal class group of $\mathbb{Q} (\sqrt{D})$ contains an element of order $g$ is bounded below by $X^{\frac{1}{2} + \epsilon(g) - \epsilon}$, where the exponent is the same as in Soundararajan's theorem. Combining this with a theorem of Frey, we give a lower bound for the number of quadratic twists of certain elliptic curves with $p$-Selmer group of rank at least $2$, where $p \in \{3,5,7\}$.
\end{abstract}

\section{Introduction}
Questions about the structure of ideal class groups are notorious and arise frequently in algebraic and analytic number theory. The Cohen-Lenstra heuristics \cite{CL} predict that the probability that the ideal class group of an imaginary quadratic field has non-trivial $\ell$ torsion for a given prime $\ell$ is
\begin{equation}\label{eq:CL}
1 - \prod_{n = 1}^{\infty} \left( 1 - \frac{1}{\ell^{n}} \right) =  \frac{1}{\ell} + \frac{1}{\ell^2} - \frac{1}{\ell^5} \cdots
\end{equation}

There is a large disparity between what was conjectured by Cohen and Lenstra and what is known. The state of the art in the direction of equation (\ref{eq:CL}) is a famous theorem of Soundararajan \cite{Sound}, who used sieving techniques to count ideals of any given order. Given an even integer $g \ge 4$, Soundararajan showed the following:
$$
\# \{ -X < D < 0 : Cl( D) \mbox{ contains an element of order $g$} \} \gg X^{\frac{1}{2} + \epsilon(g) - \epsilon},
$$
where $D$ is square free, $Cl( D)$ is the ideal class group for $\mathbb{Q} (\sqrt{D})$ and $\epsilon (g) >0$ is given by
\begin{equation}\label{eq:exponent}
\epsilon(g) = \begin{cases} 
     \frac{2}{g}  & g \equiv 0 \pmod{4} \\ 
      \frac{3}{g+2}  & g \equiv 2 \pmod{4}. \\ 
   \end{cases}
\end{equation}

Note that this theorem gives an estimate for odd order torsion as well, since a group with an element of order $2k$ must also have an element of order $k$. This theorem gives the strongest estimate for class number divisibility for most $g$. Gauss' genus theory shows that $2$-torsion subgroup of $Cl(D)$ has size $2^{w(D)}$, where $w(D)$ is the number of distinct prime factors of $D$. Furthermore, \cite{Sound} gives a lower of estimate of $X / \log{X}$ for the number of imaginary quadratic fields with an element of order 4 in their class group. For 3-torsion, Heath-Brown \cite{HB} proved the lower bound of $X^{\frac{9}{10} - \epsilon}$. 


In her PhD thesis \cite{Beckwith}, the author studied the indivisibility of class groups for imaginary quadratic fields satisfying arbitrary local conditions, refining and quantifying a recent theorem of Wiles \cite{Wiles}. The theorem of that paper was analogous to the work of Horie and Nakagawa \cite{HN}, who extended a famous theorem of Davenport and Heilbronn \cite{DH} to count imaginary quadratic fields satisfying arbitrary local conditions at primes and with trivial 3-torsion in their class group. The theorem of Horie and Nakagawa has been used to prove some of the strongest known theorems on rank one quadratic twists of elliptic curves whose Mordell-Weil group contains a 3-torsion point (see \cite{Byeon}, \cite{Vatsal}). 

In this note, instead of indivisibility we consider divisibility of class numbers, and specifically in analogy with the author's thesis work we refine Soundararajan's \cite{Sound} famous theorem. 

Throughout the paper, we assume that $A,B$ are fixed integers such that $B > 1$ and $\gcd(A,B)$ is square-free. We let $g$ be an even integer and $N_g(X;A,B)$ be the number of $0 \le D \le X$ such that 
\begin{enumerate}
\item  $Cl(-D)$ has an element of order $g$, 
\item $D \equiv A \pmod{B}$,
\item $D$ is square-free.
\end{enumerate}

When $4 \nmid g$, we are able to use small modifications of Soundararajan's argument to show that $N_g(X;A,B)$ satisfies the same lower bound as in \cite{Sound}. However, when $g \equiv 0 \pmod{4}$, the Diophantine equations in \cite{Sound} do not always have solutions modulo $B$ for square-free discriminants in arbitrarily given residue classes. For this reason, we say that a pair $(A,B)$ is \emph{special} for $g$ if the relation $2m^{g/2} - t^2 \equiv A \pmod{B}$ has a solution with $\gcd(t,B) = 1$ and $\gcd(m,2t) = 1$. 


\begin{theorem}\label{thm:divisibility}
Let $g \ge 4$ be even and let $A,B$ be fixed integers, $B > 1$. Assume that $\gcd(A,B)$ is square-free. If $g \equiv 2 \pmod{4}$ or if $(A,B)$ is special for $g$, then for all $\epsilon > 0$, we have
$$
N_g(X;A,B) \gg  X^{\frac{1}{2} + \epsilon(g) - \epsilon}. 
$$

Suppose there exist integers $x,y$ such that $0 < 2x^2 - y^2 \equiv A \pmod{B}$ is square-free and $\gcd(x,2y) = 1$. Then
$$
N_4(X;A,B) \gg \frac{X}{\log{X}}.
$$
		
\end{theorem}

Going back to the 1980s, the first major theorems on the Birch and Swinnerton-Dyer conjecture involved the study of quadratic twists of elliptic curves and imaginary quadratic fields. For example, the Gross-Zagier formula \cite{GrZa} gives a relationship between the height of Heegner points and the vanishing of $L$-functions, and when the analytic rank is one, this allows you to construct a rational Heegner point of infinite order from a suitable discriminant. The work of Kolvagin \cite{Ko1} shows that if the analytic rank of an elliptic curve $E$ over $\mathbb{Q}$ is 0, and there is a Heegner point of infinite order in $E(K)$ for an imaginary quadratic field $K$, then both the Mordell-Weil and Shafarevich-Tate groups of $E(\mathbb{Q})$ are finite.  Combining this with the Gross-Zagier formula and an analytic number theory result of Bump, Friedberg, and Hoffstein on the existence of twists whose $L$-functions have a zero of order 1 at the central value, one has that the weak Birch and Swinnerton-Dyer conjecture holds when the analytic rank is 0. 

It is natural to wonder what else is known about the structure of Mordell-Weil groups over a family of quadratic twists. Goldfeld famously conjectured that half of the quadratic twists have a Mordell-Weil group of rank 1, and half have rank 0, but this has not been shown in any cases.

Regarding larger rank, Gouv\^{e}a and Mazur \cite{GM} used sieving techniques to show that the number of twists with even analytic rank $\ge 2$ with twisting discriminant less than $X$ in absolute value is greater than $X^{\frac{1}{2} - \epsilon}$ for sufficiently large $X$. Under the assumption of the Parity Conjecture (which says that the algebraic and analytic ranks of an elliptic curve over $\mathbb{Q}$ are equal modulo 2), they showed the same lower bound for the number of twists of $E$ with rank at least 2 in their Mordell-Weil group. Stewart and Top \cite{ST} gave the unconditional lower bound of $O(X^{\frac{1}{7}} / (\log X)^2)$ for the number of twists of Mordell-Weil rank $\ge 2$.

Here, we give an estimate for the number of twists of certain elliptic curves with rank at least $2$ in their $p$-Selmer group, where $p \in \{ 3,5,7 \}$ is fixed. We apply Theorem \ref{thm:divisibility} and a theorem of Frey \cite{Frey} on twists of elliptic curves with $p$-torsion to show that under certain conditions we can produce nontrivial elements in Selmer groups in families of elliptic curves given by quadratic twists. To this end, it was required to refine Soundararajan's thesis similar to the author's refinement of Wiles' theorem. We combine this with important work of the Dokchitsers \cite{DD} establishing that the analytic rank of $E / \mathbb{Q}$ shares the same parity as the parity of the $p$-Selmer rank of $E$. 

Let $M_{E,p}(X)$ be the number of square-free $D$ with $-X < D < 0$ such that the $p$-Selmer group of the quadratic twist $E_D$ has an even rank of at least 2.

To state our result precisely, for any elliptic curve $E / \mathbb{Q}$ with a torsion point of odd prime order $p$, we let $\widetilde{S_E}$ denote the set of odd primes $q$ dividing the conductor of $E$ such that $q \equiv -1 \pmod{p}$ and $v_q(\Delta_E) \neq 0 \pmod{p}$ and $v_q(j(E)) < 0$. 

\begin{corollary}\label{thm:cor1}
Let $E/\mathbb{Q}$ be an elliptic curve. Suppose that $E$ contains a point $P$ of order $p \in \{3,5,7 \}$. Furthermore, assume that either $E$ is given by $y^2 = x^3 + 1$ or that $P$ is not contained in the kernel of the reduction modulo $p$.  Assume that $\widetilde{S_{E}} = \emptyset$. 

We assume that one of the following three conditions is true.
\begin{enumerate}
\item $N(E)$ is odd, $2 \nmid v_{p}(N(E))$ and $v_{p} (j(E)) \ge 0$ and $p | N(E)$.

\item $N(E)$ is odd and $(-1)^{\alpha(E)} = - sign ( E)$, where $\alpha(E)$ is the number of prime divisors $q$ of the square-free part of $N(E)$ such that $q \neq p$ and either $v_q(j(E)) \ge 0$ or $E (\mathbb{Q}_q)$ is a Tate curve. 

\item There exists a square-free $d > 0$ such that $N(E_d)$ has at least one prime divisor that does not divide $N(E)$, and $d \equiv 3 \pmod{4}$ if $2 | \gcd(N(E),N(E_d))$.   

\end{enumerate}

Then we have
$$
M_{E,p} (X) \gg X^{\frac{1}{2} + \frac{3}{2 p  +2} - \epsilon}.
$$  
\end{corollary}

We will prove Theorem \ref{thm:divisibility} in Sections 2-4. In Section 5 we will state and prove the corollary of our result and give some examples.
\section*{Acknowledgements}
The author thanks Ken Ono for suggesting this project and Professor Soundararajan for a helpful conversation. 

\section{Proof of Theorem \ref{thm:divisibility}}
We prove Theorem \ref{thm:divisibility} essentially by recapitulating Soundararajan's method \cite{Sound}. Soundararajan shows that solutions of certain Diophantine equations produce ideals of order $g$, and counts the solutions of these equations. To count discriminants in specific arithmetic progressions, we impose congruence conditions on the solutions of these equations, a strategy analogous to that of Horie and Nakagawa \cite{HN}. 

\textbf{Case 1: $g \equiv 2 \pmod{4}$. }

We use Proposition 1 of \cite{Sound}. 

\begin{proposition}[Soundararajan, Proposition 1]\label{thm:sound1}
Let $h \ge 3$ be an integer. Suppose that $D \ge 63$ is a square-free integer such that 
\begin{equation}\label{eq:diophantine}
t^2 D = m^{h} - n^2,
\end{equation} 
where $m,n,t$ are positive integers satisfying $\gcd(m,2n) = 1$ and $m^{h} < (D + 1)^2$. Then $Cl(-D)$ has an element of order $h$. If $h$ is odd and $D$ has at least two odd prime factors, then $Cl(-D)$ has an element of order $2h$. 
\end{proposition}


Let $g_1 = g/2$. We count solutions of equation (\ref{eq:diophantine}) for $h=g_1$ lying in specific intervals and satisfying certain congruence conditions. The intervals are the same as in \cite{Sound}, and are given below in terms of a parameter $T$ which is at most  $\sqrt{X} / 2^{g_1 +3}$, as well as $M := T^{2 / g_1} X^{1 / g_1} / 2$, and $N := TX^{\frac{1}{2}} / 2^{g_1 + 1}$.

To define the congruence conditions, we fix $r$ to be 1 if $\gcd(A,B) > 2$, and otherwise we let $r$ be a fixed odd prime that does not divide $B$. 
We choose $A',B'$ so that all integers in the residue class $A' \pmod{B'}$ are divisible by $r$, are not divisible 4 or by the square of any prime divisor of $Br$, and are in the residue class $A \pmod{B}$. This can be done by choosing an appropriate residue class modulo the square of each prime divisor of $Br$.

We let $m_1,n_1,t_1$ be as in Lemma \ref{thm:oddsolutions} for $a = -A'$ and $b = B'$, and we let $\mathbf{x} := (m_1,n_1,t_1)$. Additionally, we let $R_1 (X,T;D,\textbf{x})$ denote the number of triplets of integers $(m,n,t)$ satisfying the following:
\begin{enumerate}
\item $m \in (M, 2M]$, $t \in ( T, 2T]$, $n \in (N, 2N]$, 
\item $ m^{g_1} - n^2 = Dt^2$,  
\item $\gcd(m,2n) = 1$, 
\item $(m,n,t) \equiv \mathbf{x} \pmod{B'}$. 
\end{enumerate}

We will prove the following estimate in the next section. This is the same as the estimate for the analogous sum in \cite{Sound}.

\begin{proposition}\label{thm:oddestimate}
Let $X,T, \mathbf{x}$ and $ R_1 (X, T; D, \mathbf{x})$ be as in the preceding discussion.
\begin{align*}
\sum_{\mbox{square-free } D \le X} R_1 (X, T; D, \mathbf{x}) &\asymp \frac{MN}{T} + o(MT^{\frac{2}{3}} X^{\frac{1}{3}}) &= X^{\frac{1}{2} + \frac{1}{g_1}} T^{\frac{2}{g_1}} \\
+ o (X^{\frac{1}{3} + \frac{1}{g_1}} T^{\frac{2}{3} + \frac{2}{g_1}}).
\end{align*}
\end{proposition}

Theorem \ref{thm:divisibility} follows from Proposition \ref{thm:oddestimate} by the argument as in \cite{Sound}. Let $S_{1,g} (X, T; \mathbf{x})$ be the number of square free $1 \le D \le X$ such that there exists a triplet $(m,n,t)$ satisfying (1) - (4) in the definition of $R_1(X,T;D,\mathbf{x})$. Because of conditions (2) and (4), for such $D$ we have $D \equiv -A' \pmod{B'}$. As a consequence of (1), we have $m^{g_1} < (D+1)^2$. By Proposition \ref{thm:sound1}, we have $N_g (X; A,B) \ge S_{1,g} (X,T;\mathbf{x})$. 


By the same combinatorial argument as in (1.5) of \cite{Sound}, in which one bounds the number of pairs of distinct triplets of solutions $(m,n,t)$ for $1 \le D \le X$, we have:

$$
\sum_{D \le X} R_1 (X, T; D,\mathbf{x}) ( R_1 (X, T; D,\mathbf{x}) - 1) \ll T^{2 + \frac{4}{g_1}} X^{\frac{2}{g_1} + \epsilon}.
$$

The Cauchy-Schwarz inequality implies
$$
S_{1,g} (X, T; A,B) \ge \left( \sum_{D \le X} R_1 (X, T; D, \mathbf{x})  \right)^2 \left( \sum_{D \le X} R_1 (X, T;D, \mathbf{x})^2   \right)^{-1}.
$$

Theorem \ref{thm:divisibility} for $g \equiv 2 \pmod{4}$ follows from setting $T = X^{\frac{ g_1 -2 }{4(g_1 + 1)}}$. 
 
\ \\ \

\textbf{Case 2: $g \equiv 0 \pmod{4}$.}

We use the following result of \cite{Sound}. 

\begin{proposition}[Soundararajan \cite{Sound}, Proposition 1']\label{thm:evenequation}
Suppose that $g \ge 3$ is even and that $D = 2 m^{g/2} -  t^2 \ge 63$ is square-free, where $\gcd(m,2t) = 1$ and $m^{g/2} < D+1$. Then $Cl(-D)$ has an element of order $g$. 
\end{proposition}


As in Case 1, we choose an arithmetic progression $A' \pmod{B'}$ to ensure that the  discriminants appearing are not divisible by $p^2$ for any prime $p$ that divides $B$, and so that $(A',B')$ is also special. 

We let $m_0,t_0$ satisfy $2m_0^{g_1} - t_0^2 \equiv -A' \pmod{B'}$ with $\gcd(B',t_0) = \gcd(m_0,2t_0) = 1$. 

Then we have $N_g (X; A,B) \ge S_{0,g} (X;m_0,t_0)$, where $S_{0,g} (X; m_0, t_0)$ is the number of square-free $D \equiv -A' \pmod{B'}$ between 0 and $X$ such that there are integers $(m,t)$ satisfying
\begin{enumerate}
\item $D = 2 m^{g/2} - t^2$,
\item $(m,t) \equiv (m_0, t_0) \pmod{B'}$,
\item $ m^{g_1} < D+1$,
\item $\gcd(m, 2t ) = 1$. 
\end{enumerate}
We define $R_0 (D; m_0,t_0)$ to be the number of pairs $(m,t)$ satisfying (1) - (4). By Proposition 2 of \cite{Sound}, we have $R_0 (D; m_0, t_0) \le \tau(D) /2 \ll D^{\epsilon}$.

To conclude Theorem \ref{thm:divisibility} for Case 2, it is sufficient to show the following.
\begin{proposition}\label{thm:evenestimate}
Let $A,B,$ and $R_0(D;t_0)$ be as above. 
\begin{equation}\label{eq:R0}
\sum_{\mbox{square-free }D \le X} R_0 (D; m_0,t_0) \gg X^{\frac{1}{2} + \frac{1}{g_1}}.
\end{equation}
\end{proposition}

We prove Proposition \ref{thm:evenestimate} in the next section.

\textbf{Case 3: $g = 4$ }

The claim follows from Propositions 1' and 2 of \cite{Sound}. Specifically, Proposition 1' says that if $D = 2x^2 - y^2$ is square-free and $\gcd(x,2y) = 1$, then $Cl(-D)$ has an element of order 4. Proposition 2 \cite{Sound} says that for a given square-free $D$, the number of such pairs $(x,y)$ is $\tau(D) / 2$ if $D \equiv 1 \pmod{8}$ and $D$ is divisible only by primes which are $\pm 1 \pmod{8}$, and that the number of such pairs is 0 otherwise. We see that for any prime $p > \gcd(A,B)$ for which $p \cdot \gcd(A,B) \equiv 1 \pmod{8}$ and $p  \equiv -A/\gcd(A,B) \pmod{B/\gcd(A,B)}$, we have that $Cl(-p \cdot \gcd(A,B))$ has an element of order 4, and the result follows from Dirichlet's theorem on primes in arithmetic progressions. 

\qed 
\section{Counting Diophantine solutions} 

To finish Theorem \ref{thm:divisibility}, it remains to prove Propositions \ref{thm:oddestimate}and \ref{thm:evenestimate}, which are analogous to relations (1.4) and (1.2) of \cite{Sound}, respectively.

\textbf{Proof of Proposition \ref{thm:oddestimate}}

To estimate the desired sum, we will consider three sets of integer triplets. We let $N_1 = N_1(X,T;\mathbf{x})$ be the number of 3-tuples $(m,n,t)$ satisfying conditions (1)-(4) of Proposition \ref{thm:oddestimate} for some $D$ that is indivisible by $p^2$ for all primes $p \le \log X$. We let $N_2 := N_2(X,T;\mathbf{x})$ be the number of tuples also satisfying (1)-(4) such that $D$ is divisible by $p^2$ for some $\log X \le p \le Z := X^{\frac{1}{3}} T^{- \frac{1}{3}} (\log X)^{\frac{2}{3}}$, and similarly, we let $N_3 := N_3(X,T;\textbf{x})$ be the number of tuples with $\frac{m^{g_1} - n^2}{t^2}$ divisible by $p^2$ for some $p > Z$. 

The same estimates as in \cite{Sound} hold for $N_2$ and $N_3$ by the same argument, specifically,
$$
N_2 \ll \frac{MN}{T \log X} + o(MX^{\frac{1}{3}}) \quad \mbox{and} \quad N_3 \ll \frac{X}{Z^2} M \log{X} = o (M X^{\frac{1}{3}} T^{\frac{2}{3}} ).
$$

We have $0 \le N_1 - \sum_{D \le X} R_1(X,T; \textbf{x}) \le N_2 + N_3$. The rest of the section is devoted to estimating $N_1$.


For fixed $(m,t) \equiv (m_1,t_1) \pmod{B'}$ such that $m \in (M,2M]$ is odd and $t \in (T,2T]$, let $N_1 (m,t)$ be the number of $n \in (N, 2N]$ such that $n \equiv n_1 \pmod{B'}$ and $\frac{m^{g_1} - n^2}{t^2}$ is indivisible by $p^2$ for all primes $p \le \log X$ with $p \nmid B'$. 

Let $P$ be given by $P : =  \prod_{\substack{ 2 \le p \le \log{x} \\ p \nmid B' } } p$. 

\begin{align*}
N_1 (m,t) &= \sum_{\substack{ N \le n \le 2N, (n,m) = 1 \\ n \equiv n_1\pmod{B'} \\ n^2 \equiv m^{g_1} \pmod{t^2}  }}   \left(\sum_{\ell^2 | (\frac{m^{g_1} - n^2}{t^2}, P^2)} \mu (\ell) \right) \\
&= \sum_{\ell | P, (\ell,m) = 1} \mu(\ell) \sum_{\substack{ N \le n \le 2N \\ n \equiv n_1\pmod{B'} \\ n^2 \equiv m^g_1 \pmod{t^2 \ell^2} }} 1
\end{align*}

We evaluate the inner sum by decomposing $[N,2N]$ into subintervals of length $B' t^2 \ell^2$. 
$$
\sum_{\substack{ N \le n \le 2N \\ n \equiv n_0 \pmod{B'} \\ n^2 \equiv m^{g_1} \pmod{t^2 \ell^2}  }} 1 
= \frac{N}{B't^2 \ell^2} \left( \sum_{ \substack{n \equiv n_1 \pmod{B'} \\ n^2 \equiv m^{g_1} \pmod{t^2 \ell^2}  \\ n \pmod{B' t^2 \ell^2}}} 1 \right)  + O(\varphi_m(t^2 \ell^2 B')) 
$$

We let $\varphi_m(\ell)$ be the number of solutions modulo $\ell$ to the congruence $n^2 \equiv m^{g_1} \pmod{\ell}$. As $n$ varies over the residue classes modulo $t^2 \ell^2 B'$ such that $n \equiv n_0 \pmod{B'}$, we have that $n$ will also vary over each residue class modulo $\ell^2 t^2$. Therefore, the value of the summation on the right hand side is $\varphi_m(t^2 \ell^2)$.

$$
\sum_{\substack{ N \le n \le 2N \\ n \equiv n_1 \pmod{B'} \\ n^2 \equiv m^{g_1} \pmod{t^2}  }} 1 
= \frac{N}{B't^2 \ell^2} \varphi_m(t^2 \ell^2)  + O(\varphi_m(t^2 \ell^2 B')). 
$$ 

Next we evaluate $N_1 (m,t)$ using the properties of $\varphi_m(n)$ as a function of $n$. In particular, $\varphi_m$ is a multiplicative function such that for any prime $p \nmid m$, we have $\varphi_m(p^n) = \varphi_m(p) = 1 + \left(\frac{m}{p} \right)$. 


\begin{align*}
N_1(m,t) &= \frac{N}{B't^2} \sum_{\ell | P, (\ell,m)=1} \mu(\ell) \cdot  \frac{1}{\ell^2} \varphi_m ( t^2 \ell^2) + O(X^{\epsilon} \tau(P)) \\
&= \frac{N}{B't^2} \varphi_m(t^2) \sum_{\ell | P, (\ell,m)=1} \mu(\ell) \cdot  \frac{1}{\ell^2} \varphi_m ( \ell / (t, \ell)) + O(X^{\epsilon} \tau(P)) \\
&= \frac{N}{B't^2} \varphi_m(t^2) \prod_{p|P, p \nmid m} \left( 1-  \frac{1}{p^2} \varphi_m ( p / (t, p)) \right) + O(X^{\epsilon}) \\
&\asymp \frac{N}{T^2} \varphi_m (t^2) + O(X^{\epsilon}),
\end{align*}
where the implied constants in the last expression are of course independent of $X$,$m,t$.  


Soundararajan shows (Lemmas 1 and 2 \cite{Sound}) that the average value of $\varphi_{m} (t^2)$ is 1 as $t$ and $m$ vary over the intervals and residue classes appearing in the definition of $N_1$. 

\begin{lemma}\label{thm:averagemodm}
Assume $g_1$ is odd, and $a,b,c$ are integers such that $(ab,c) = 1$.    
\begin{align*}
\sum_{\substack{M < m \le 2M \\ m \equiv a \pmod{c}}} \sum_{\substack{T < t \le 2T \\ \gcd(t,m) = 1 \\ t \equiv b \pmod{c}}}\varphi_m (t^2) &= \sum_{\substack{M < m \le 2M \\ m \equiv a \pmod{c}}} \sum_{\substack{T < t \le 2T \\ \gcd(t,m) = 1 \\ t \equiv b \pmod{c}}} 1 + O(TM^{\frac{5}{8}} (\log X)^3)   \\
&\asymp MT + O(TM^{\frac{5}{8}} (\log X)^3).
\end{align*}
\end{lemma}

The proof, which uses the identity $\varphi_m(t^2) = \sum_{d | t} \mu^2(d) \left( \frac{m}{d} \right)$ and the P\'{o}lya-Vinogradov inequality, is identical to the proof of Lemma 2 of \cite{Sound}, so we omit it for the sake of brevity. 

Now we have 
\begin{align*}
N_1(X,T;\textbf{x}) &= \sum_{\substack{T < t \le 2T \\ M < m \le 2M \\ (t,m) \equiv (t_1,m_1) \pmod{B'}}} N_1(m,t) \\
&\asymp \frac{N}{T^2} \sum_{\substack{T < t \le 2T \\ M < m \le 2M \\ (t,m) \equiv (t_1,m_1) \pmod{B'}}} \varphi_m (t^2) + O(X^{\epsilon} MT) \\
&\asymp \frac{NM}{T}  + o( M^{\frac{1}{3}} T^{\frac{1}{3}})
\end{align*}

\textbf{Proof of Proposition \ref{thm:evenestimate}}
The strategy is a little bit simpler than the previous proof. We let $M_1$ be the number of pairs $(m,t)$ satisfying (1)-(4) of Proposition \ref{thm:evenestimate} and also satisfying the property that for all $p \le \log{X}$, $p^2 \nmid 2 m^{g_1} - t^2$. We let $M_2$ be the number of pairs $(m,t)$ satisfying (1)-(5) and also satisfying the property that $2m^{g_1} - t^2$ is divisible by $p^2$ for some $p > \log X$. The desired sum is $M_1 - M_2$.

To estimate $M_1$, we fix $m \equiv m_0 \pmod{B'}$ such that $2 \nmid m$. We denote by $M_1(m;t_0)$ the number of values $t \le X$ such that $(m,t)$ satisfy the requirements of $M_1$. 

We have
\begin{align*}
M_1 (m;t_0) &= \sum_{\substack{\sqrt{X}/2^{g_1 + 1} \le t \le \sqrt{X}/2^{g_1} \\ (t,m) = 1, t \equiv t_0 \pmod{B'}}} \sum_{\ell^2 | \gcd(2m^{g_1} - t^2,P^2)} \mu(\ell) \\
&= \sum_{\ell | P, 1=(\ell,m)} \mu(\ell) \sum_{\substack{\sqrt{X}/2^{g_1 + 1} \le t \le \sqrt{X}/2^{g_1} \\ t \equiv t_0 \pmod{B'} \\ 2 m^{g_1} \equiv  t^2 \pmod{\ell^2}}} 1 \\
&= \sum_{\ell | P, 1=(\ell,m)} \mu (\ell)  \left( \frac{\sqrt{X}}{2^{g_1 + 1} \ell^2} \varphi_2' (\ell^2) + O( \varphi_{2}'(B'\ell^2)) \right) \\
&= \frac{\sqrt{X}}{2^{g_1 + 1}} \prod_{p | P, p\nmid m} (1 -   \frac{1}{p^2} \varphi_2' (p^2))  + O(X^{\epsilon}) \\
&\asymp \sqrt{X}
\end{align*}
 
Here $\varphi_m'(t)$ is the number of $n \pmod {t}$ such that $n^2 \equiv m \pmod{t}$. 

So we have 
\begin{align*}
M_1 &= \sum_{\substack{m \equiv m_1 \pmod{B'}, 2 \nmid m \\ X^{1/g_1} / 4 \le m \le X^{\frac{1}{g_1}}}} M_1(m;t_0) \\
&\asymp X^{\frac{1}{2} + \frac{1}{g_1}}.
\end{align*}

By the same argument as in \cite{Sound} for $M_2$, we have $M_2 \ll \frac{X^{\frac{1}{2} + \frac{1}{g_1}}}{\log X}$. This concludes the proof of Proposition \ref{thm:evenestimate}. 


\section{Diophantine equations modulo $m$}

The two elementary lemmas given here establish the existence of solutions lying in the necessary residue classes of the equations appearing in the propositions of the previous section.

\begin{lemma}\label{thm:oddsolutions}
Let $a,b$ be integers. If $g_1$ is odd, or if $g_1$ is even $b$ is odd, then there exist $n_1,t_1,m_1$ such that $(b,t_1) = 1 = (m,2n) = (m,b)$, and $m_1^{g_1} - n_1^2 \equiv t_1^2 a \pmod{b}$. 
\end{lemma}


From Hensel's lemma and the Chinese Remainder Theorem, it is sufficient to consider the case that $p$ is prime. For prime $b \ge 7$, we use the fact given below as Lemma \ref{thm:fewsquares}. We let $t_1 \equiv v^{-1} \pmod{b}$ and $n_1 \equiv u t_1 \pmod{b}$ and $m_1 \equiv 1 \pmod{b}$. 

For $b = 2,3,5$ we check each case of $a \pmod{b}$ and find that an appropriate solution exists. 

\qed

\begin{lemma}\label{thm:fewsquares}
Suppose $A$ is an integer, and $p \ge 7$ is prime. Then there exists $u,v$ such that $A + u^2 \equiv v^2 \pmod{p}$  and $\gcd(uv,p)=1$.

\end{lemma}
Proof: 
If $A \equiv 0 \pmod{p}$, then we take $u = v=1$. Otherwise, we can find explicit choices easily using the factorization $v^2 - u^2 = (v - u)(v + u)$. For example, $v \equiv 2^{-1} (A + 1)$ and $u \equiv A - 2^{-1} - 2^{-1} A$ is an appropriate solution unless $p$ divides one of these choices, in which case $A \equiv \pm 1 \pmod{p}$, in which case $v \equiv 1 + 2^{-2} A$ and $u \equiv -1 + 2^{-2} A$ is an appropriate choice. 

\qed

\section{Selmer groups}

\subsection{Frey's Theorem}

For any elliptic curve $E / \mathbb{Q}$, we let $\Delta_E$, $N(E)$, and $j(E)$ denote the discriminant, conductor, and $j$-invariant for $E$, respectively. We let $Cl(D)_{p}$ denote the $p$-torsion subgroup of $Cl(D)$. 

Frey \cite{Frey} proved a double divisibility relationship between $ Cl(D)_p$ and the size of the $p$-Selmer group of $E_D$, denoted $Sel_{p} (E_d, \mathbb{Q}) $. A simplified version of his theorem is stated below, providing conditions under which $\# Cl(D)_p$ divides $\# Sel_p(E_d,\mathbb{Q})$.

\begin{theorem}[Frey]\label{thm:Frey}
Let $E$ be an elliptic curve defined over $\mathbb{Q}$ with a rational point $P$ of prime order $p \in \{3,5,7 \}$ over $\mathbb{Q}$. Assume that either $E$ is given by the equation $y^2 = x^3 + 1$ (hence $p =3$) or that $P$ is not contained in the kernel of reduction modulo $p$. 

Also assume that $\widetilde{S_E} = \emptyset$. 

Let $d$ be a square-free integer relatively prime to $pN(E)$ such that
\begin{enumerate}
\item If $2 | N_E$ then $d \equiv 3 \pmod{4}$
\item If $q$ is a prime divisor of $N(E)$ other than $2$ or $p$, then $\left( \frac{d}{q} \right) = -1$ if $E$ is a Tate curve over $\mathbb{Q}_q$ or $v_q(j_E) \ge 0$, and $\left( \frac{d}{q} \right) = 1$ otherwise. 
\item If $v_p(j_E) < 0 $ then $\left( \frac{d}{q} \right) = -1$. 
\end{enumerate}

Then we have 
$$
\# Cl(d) | \# Sel_p(E_d, \mathbb{Q}) .
$$ 
\end{theorem}

\subsection{Proof of Corollary \ref{thm:cor1}}
We choose $A$, $B$ in the following way so that each square free $D \equiv A \pmod{B}$ satisfies the conditions of Frey's theorem for $E$. For each prime $q \neq \ell$ dividing the conductor of $N({E})$, if $E$ is a Tate-curve over $\mathbb{Q}_p$ or $\nu_q(j_E) \ge 0$, then we let $q | B \nmid q^2$ and $\left( \frac{A}{q} \right) = -1$. For all other $q \neq \ell$ dividing the conductor of $E$, we let $q | B \nmid q^2$ and $\left( \frac{A}{q} \right) = 1$.  If we are in Case 1, we choose $A$ so that $\left( \frac{A }{p} \right) = \sgn(E) (-1)^{\alpha(N(E))}$. If we are in Case 2, we choose $A \pmod{B}$ so that $\left( \frac{D}{p} \right) \neq =0$.  In Case 3, we also require that $D$ is divisible by $d$, and that $\gcd(D/d, N(E_d)) = 1$. Furthermore, we choose conditions modulo the prime divisors of $N(E_d)$ to ensure that we have
$$
\left( \frac{D/d}{N(E_d)} \right) = - \sgn(E_d).
$$

For all $D \equiv A \pmod{B}$, the conditions of Frey's theorems are met. 

It also turns out that for such $D$, we have $\sgn(E_D)  = 1$. This is because we have, for fundamental discriminants $D$ with $\gcd(D,N(E) ) = 1$, we have that
$$
\sgn(E_D) = \sgn(E) \chi_D (-N(E)),
$$
where $\chi_D$ is the character associated to $\mathbb{Q}(\sqrt{D})$. 

In the first case, the right hand side is $- \sgn(E) (-1)^{\alpha(N(E))} \left( \frac{D}{p} \right)$. By our choice of $\left( \frac{D}{p} \right)$, this is 1. In the second case, the right hand side is simply $- \sgn(E) (-1)^{\alpha(N(E))}  =1$.

In the third case, the conditions modulo 4 ensure that $D / d \equiv 1 \pmod{4}$ if $2 | N(E_d)$. We then compute 
$$
\sgn(E_D) = \sgn(E_d) \chi_{D/d} (- N(E_d)) = 1.
$$

%

For twists $E_D$ with sign 1, the $p$-rank of $Sel_p(E_{D},\mathbb{Q})$ is even. By Theorem 1.4 of \cite{DD}, the $p^{\infty}$ Selmer-rank  of $E_D$  is even. No nontrivial twists of $E$ contain $p$-torsion. It follows that for $D \neq 1$, the $p$-rank of $Sel_{p}(E_D,\mathbb{Q})$ is equal to the ${p}^{\infty}$-Selmer rank of $E / \mathbb{Q}$, and therefore the $p$-rank of $Sel_p(E_D,\mathbb{Q})$ is even. 
From the above discussion, for $D \equiv A \pmod{B}$, the $p$-rank of the Selmer group of $E_D$ is even and is at least 1 by Theorem \ref{thm:Frey}. 

The corollary now follows from Theorem \ref{thm:divisibility}. 

\subsection{Examples}
We consider some elliptic curves which satisfy the conditions of Corollary \ref{thm:cor1}. 

\textbf{Example 1:} Cremona 27a4, which is defined by $E:y^2+y = x^3 - 30x +63$ has the 3-torsion point $P=(3:0:1)$. One can verify that $E$ satisfies the conditions of Corollary \ref{thm:cor1} (1) by computing the $j$-invariant or consulting lmfdb.org. It follows that for any $\epsilon > 0$, we have
$$
M_{E,3}(X) \gg X^{\frac{7}{8} - \epsilon} . 
$$

\textbf{Example 2:} Cremona 175a2 is defined by $E: y^2 + y = x^3 - x^2 - 148x + 748$. The 5-torsion point $(-13:17:1)$ lies outside the kernel of reduction modulo 5, and the other conditions of Corollary \ref{thm:cor1} (1) can be verified. We conclude that we have
$$
M_{E,5}(X) \gg X^{\frac{3}{4} - \epsilon}
$$
for any $\epsilon > 0$.

\textbf{Example 3:} We consider Cremona 574i1, which has the $7$-torsion point $(2: 11 :1)$. The twist $E_3$ has conductor $574\cdot 3^2 \cdot 2^3$. By Corollary \ref{thm:cor1} (3) with $d =3$, we have

$$
M_{E,7}(X) \gg X^{\frac{11}{16} - \epsilon}
$$
for any $\epsilon > 0$.

\end{document}